\documentclass{article}
\usepackage{amsmath,amssymb}

\usepackage{tikz}
\usetikzlibrary{calc,decorations.markings,arrows.meta}
\tikzset{baseline={([yshift=-.5ex]current bounding box.center)}}
\tikzstyle directed=[blue,postaction={decorate,decoration={markings, mark= at position .55 with {\arrow{stealth}}}}]
\tikzstyle{v} = [draw,black,fill,circle, inner sep=1pt]
\tikzstyle{->}=[-{Stealth[]}]
\tikzstyle{<->}=[{Stealth[]}-{Stealth[]}]
\tikzstyle{<-}=[{Stealth[]}-]

\usepackage{graphicx}

\def\C{{\mathbb C}}
\def\S{{\mathbb S}}

\def\cA{{\mathcal A}}
\def\cB{{\mathcal B}}

\def\cG{{\mathcal G}}

\def\gl{{\mathfrak{gl}}}
\def\sl{{\mathfrak{sl}}}
\def\so{{\mathfrak{so}}}

\newtheorem{theorem}{Theorem}[section]
\newtheorem{lemma}[theorem]{Lemma}
\newtheorem{proposition}[theorem]{Proposition}
\newtheorem{remark}[theorem]{Remark}
\newtheorem{corollary}[theorem]{Corollary}
\newtheorem{example}[theorem]{Example}
\newtheorem{definition}[theorem]{Definition}

\title{The universal  $\gl$-weight system and the chromatic polynomial}

\author{M.~Kazarian\thanks{Higher School of Economics, Skolkovo Institute of Science and Technology
}, N.~Kodaneva \thanks{Higher School of Economics}, 
S.~Lando
\thanks{Higher School of Economics, Skolkovo Institute of Science and Technology
}}

\begin{document}

\maketitle
\date{}

\newpage

\begin{tabular}{rl}

$\gl(N)$ & Lie algebra of all $N\times N$ matrices \\
$E_{i,j}\in\gl(N)$ & matrix unit, $1\le i,j\le N$\\
$w_\gl$ & weight system associated to the series $\gl$ of Lie algebras\\
$U\gl(\infty)$ & universal enveloping algebra of the limit Lie algebra $\gl(\infty)$\\
$ZU\gl(N)$ & center of $U\gl(N)$\\
$C_1,C_2,\dots$ & Casimir elements generating the center of~$U\gl(\infty)$\\

$\alpha\in\S_m$ & permutation, which determines the hyperedges of the hypermap\\
$\S_m$ & group of permutations of $m$ elements\\
$G(\alpha)$ & the digraph of a permutation $\alpha$\\
$f(\alpha)$ & number of faces of~$\alpha$\\
$c(\alpha)$ & number of cycles in~$\alpha$\\
$V(\alpha)$ & set of cycles of $\alpha$\\
$v_1, v_2, \ldots\in V(\alpha)$  & cycles of $\alpha$\\
$a(\alpha)$ & number of ascents of $\alpha$\\
$\alpha|_U$ & restriction of~$\alpha$ to $U\subset\{1,2,\ldots m\}$\\
$\gamma(\alpha)$ & the intersection graph of~$\alpha$\\

$\sigma=(1,2,\dots,m)\in\S_m$ & standard long cycle, the hypervertex of a hypermap\\
$\varphi\in\S_m$ & permutation, which determines the hypefaces of the hypermap\\

$\cA_m$ & the vector space spanned by equivalence classes of permutations of $m$ elements\\
$[\alpha]$ & equivalence class of $\alpha$\\
$\cA$ & the rotational Hopf algebra of permutations\\
$\alpha_1 \# \alpha_2$ & concatenation product of permutations $\alpha_1$ and $\alpha_2$\\
$\mu(\alpha)$ & coproduct of $\alpha$\\
$\bar\alpha$ & canonical form of $\alpha$\\
$P$ & subspace of primitive elements of a Hopf algebra\\
$\pi$ & projection to the subspace of primitive elements\\

$\cA_m^+$ & the vector space spanned by equivalence classes\\& of positive permutations of $m$ elements\\
$\cA^+$ & the rotational Hopf algebra of positive permutations\\

$\cB_m$ & the vector space spanned by equivalence classes\\& of involutions of $m$ elements without fixed points\\
$\cB$ & the rotational Hopf algebra of chord diagrams\\

$\cG_n$ & the vector space spanned by graphs with~$n$ vertices\\
$\cG$ & the Hopf algebra of graphs\\

$G$ & simple graph\\
$K_n$ & complete graph with $n$ vertices\\
$V(G)$ & the set of vertices of~$G$\\
$G|_U$  & the subgraph of $G$ induced by~$U\subset V(G)$\\

$\chi_G(x)$ & chromatic polynomial of~$G$\\

$X(\alpha)(N;p_1,p_2,\dots)$ & prechromatic substitution\\
$X_0(\alpha)$ & top coefficient of $X_\alpha$\\

\end{tabular}

\newpage 

\section{Introduction}
Weight systems are functions on chord diagrams
satisfying so-called Vassiliev's $4$-term relations.
They are closely related to finite type knot invariants,
see~\cite{V90,K93}. 

Certain weight systems can be derived
from graph invariants, see a recent account in~\cite{KL22}.
One  of the first examples of such an invariant was the chromatic
polynomial~\cite{CDL1}.

Another main source of weight systems are Lie algebras,
the construction due to D.~Bar-Natan~\cite{BN95} 
and M.~Kontsevich~\cite{K93}.
In recent papers~\cite{ZY22,KL22}, the weight systems
associated to the Lie algebras~$\gl(N)$, $N=1,2,3,\dots$,
were unified in a
universal $\gl$-weight system, which takes values in
the ring $\C[N,C_1,C_2,C_3,\dots]$ of polynomials in
infinitely many variables. Note that this weight system
is associated to the HOMFLYPT polynomial, which is
an important and powerful knot invariant.
The  unification has been achieved by extending the
$\gl(N)$-weight systems from chord diagrams, which can be considered as involutions without fixed  points modulo
cyclic shifts, to arbitrary permutations\footnote{The authors are grateful
to M.~Karev for  pointing out  that  the extension can be considered as a special case
of  a  construction in~\cite{DK07}. }.
The main goal of the extension was to produce an efficient 
way to compute explicitly the values of the 
$\gl(N)$-weight systems, and M.~Kazarian suggested 
a recurrence relation for such a computation, which  works
for permutations rather than just for chord diagrams.

A natural question then arises, namely, which  already
known weight systems can be obtained from the universal
$\gl$ weight system.  In addition to understanding the internal
relationship between weight systems, knowing that a
given weight system can be induced from the $\gl$-weight
system  would immediately lead to extending the former
to arbitrary permutations.

To each chord diagram, one can associate a graph, called the intersection
graph of the chord diagram.
Certain weight systems are completely determined by the intersection
graphs; this is true, for example, for the $\sl(2)$-
and for the $\gl(1|1)$-weight systems~\cite{CL07}.
This is true no longer for more complicated Lie algebras, say, for $\sl(3)$.
In general, the relationship between Lie algebra weight systems and
polynomial graph invariants looks rather mysterious.

In  our recent  paper~\cite{KoL23},
we suggested a substitution that makes the  universal  $\gl$ weight
system into the interlace polynomial. The present paper
is devoted  to  investigation  of the relationship between 
the  universal $\gl$ weight system and  the chromatic polynomial.

In this case, the situation is less straightforward. Namely, 
it is easy  to show  that there  does  not  exist  a  substitution
transforming  the universal  $\gl$ weight system  into  the chromatic
polynomial. However, we suggest a substitution, from which the chromatic polynomial emerges  as
the \emph{leading  term}  of the result,  as  $N$,  the rank of the Lie  algebras
$\gl(N)$, tends  to  infinity.  This substitution  provides  an interesting
function  on  permutations, which takes  values in  the ring of  polynomials
in  two  variables. We start to  study properties of  this  function.
In particular, we prove that the leading term in~$N$
coincides with the chromatic polynomial of the intersection graph of the permutation
not only for chord diagrams, but for a wider class of permutations,
which we call positive permutations.

The existence of a relationship between the chromatic polynomial weight system
and the series $\gl(N)$ of Lie algebras was established earlier in~\cite{L00}.
Our construction seems to be much simpler and more straightforward than the one in~\cite{L00}.

The results of the present paper lead
to numerous new questions.
For example, it would be interesting to know whether
the chromatic polynomial can be naturally extended to 
hypermaps with arbitrarily  many hypervertices  (a
permutation considered modulo  cyclic shifts can be
treated as a hypermap with a single vertex).
Also, there  are universal weight systems 
for permutations constructed
through other series of Lie algebras and Lie 
superalgebras~\cite{KY23}; their specializations  
also can lead to already known, as well as new,
weight systems and their extensions to permutations.


\subsection{$4$-term relations and weight systems}

In V.~Vassiliev's theory~\cite{V90} of finite type
knot invariants, a function on chord diagrams 
with~$n$ chords is associated to every knot invariant of order not greater than~$n$. He showed that any function constructed
in this way satisfies so-called $4$-term relations.
Then in \cite{K93}, M.~Kontsevich showed that every function on chord diagrams satisfying the four-term relations defines a knot invariant of finite type.

Below, we give necessary definitions and state the properties of the objects
we require. A complete account can be found in~\cite{LZ03,CDBook12}.

\begin{figure}
    \centering
    \includegraphics[width=0.8\linewidth]{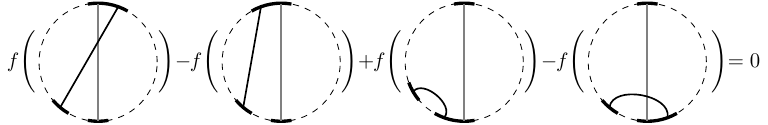}
    \caption{$4$-term relation}
    \label{fig:4t}
\end{figure}

A \emph{chord diagram} of order~$n$ is an oriented circle with $2n$~points on it, which are split into $n$ disjoint~pairs, each pair connected by a chord, considered up to orientation
preserving diffeomorphisms of the circle. 
A function~$f$ on chord diagrams is called a \emph{weight system} if it satisfies the $4$-term relation shown in Fig.~\ref{fig:4t}. The vector space spanned by chord diagrams
modulo $4$-term relations is endowed with a natural multiplication:
modulo $4$-term relations, the concatenation of two chord
diagrams is well defined.

\subsection{Graphs, embedded graphs and hypermaps}
Weight systems are closely related to invariants of graphs,
embedded ($=$ ribbon, topological) graphs and delta-matroids, see a detailed
account in~\cite{KL22}. Here we recall briefly certain notions 
that will be used in the present paper.

A chord diagram can be treated as an
orientable embedded graph (a combinatorial map)
 with a single vertex: the vertex is the outer circle, while the edges
(ribbons) are obtained by thickening the chords.

Permutations considered up to standard cyclic shift can be treated as
orientable hypermaps with a single vertex.
A \emph{hypermap} is a triple of permutations $(\sigma,\alpha,\varphi)$ in~$\S_m$ whose product is the identity permutation,
$\sigma\circ\alpha\circ\varphi={\rm id}$, considered up
to simultaneous conjugation by one and the same permutation.
A hypermap is a map if the permutation~$\alpha$ is an involution
without fixed points.
Cycles in the decomposition of $\sigma$ (resp., $\alpha$, $\varphi$) into a product
of disjoint cycles are \emph{hypervertices}
(resp., \emph{hyperedges}, \emph{hyperfaces}) of the hypermap.

Below, we will always assume that the first permutation
$\sigma$ is the standard cycle, $\sigma=(1,2,\dots,m)$,
and since $\varphi=\alpha^{-1}\circ\sigma^{-1}$ is uniquely 
reconstructed from~$\alpha$, we will speak about~$\alpha$
as about a hypermap with a single vertex.

The number
$f(\alpha)$ of (hyper)faces for a given permutation~$\alpha$
is the number of disjoint cycles in the permutation~$\varphi$. If~$\alpha$
is an involution without fixed points (an arc diagram), then the number
of its hyperfaces coincides with the  number of ordinary faces.

To each hypermap, an orientable two-dimensional surface
with boundary is naturally associated. The hyperfaces of 
the hypermap are in a natural one-to-one correspondence with
the connected components of the  boundary of the surface. 
For a single vertex, the surface is constructed by attaching
the hyperedges (which are disks corresponding to the cycles in~$\alpha$)
to the vertex, which is the disk, which corresponds to the unique cycle in~$\sigma$. The rule of gluing is prescribed by~$\alpha$.

\subsection{$\gl$-weight system}
The $\gl$-weight system~$w_\gl$ for permutations is defined in~\cite{ZY22}, see also~\cite{KL22}.
This function on permutations takes values in the algebra $\C[N][C_1,C_2,\dots]$
of polynomials in infinitely many variables $C_1,C_2,\dots$ whose coefficients depend
on the additional variable~$N=C_0$. The variables $C_i$ should be treated
as the standard \emph{Casimir elements}, which generate the center of the universal
enveloping algebra $U\gl(\infty)$ of the limit Lie algebra $\gl(\infty)$.
The definition proceeds as follows.

The value of the $\gl(N)$-weight system on a permutation~$\alpha\in \S_m$
is defined as the element
$$
w_{\gl(N)}(\alpha)=\sum_{i_1,\dots,i_m=1}^NE_{i_1,i_{\alpha(1)}}\dots E_{i_m,i_{\alpha(m)}}
$$
of the universal enveloping algebra $U\gl(N)$, where $E_{i,j}\in\gl(N)$ are the matrix units, $1\le i,j\le N$.

An arc diagram can be interpreted as an involution without fixed points on 
the set of the arc ends. 
The  $\gl$-weight system for permutations extends 
the $\gl(N)$-weight system for
chord diagrams to arbitrary permutations. 
This extension provides a way for recursive computation of the $\gl(N)$-weight
systems, as well as unifies the latter weight systems in a universal one, containing all the
specific weight systems, for a given~$N$.

In order to introduce the recurrence, we will need the notion of the digraph of a permutation.
	A permutation can be represented as an oriented graph.
The $m$ vertices of the graph correspond to the permuted elements.
They are placed on the horizontal line, and numbered from left to right
in the increasing order. The arc arrows
show the action of the permutation (so that each vertex is incident with exactly
one incoming and one outgoing arc edge).
	The digraph $G(\alpha)$ of a permutation $\alpha\in \S_m$
consists of these~$m$ vertices and $m$ oriented edges, for example:
	\[
	G(\left(1\ n+1)(2\ n+2)\cdots(n\ 2n)\right)=\begin{tikzpicture}[baseline={([yshift=-.5ex]current bounding box.center)},decoration={markings, mark= at position .55 with {\arrow{stealth}}}]
		\draw[->,thick] (-2,0)--(2,0);
		\draw[blue,postaction={decorate}] (-1.8,0) ..  controls (-1,.5) ..(.2,0);
		\draw[blue,postaction={decorate}] (-1.4,0) ..  controls (-.6,.5) ..(.6,0);
		\draw[blue,postaction={decorate}] (-.2,0) ..  controls (1,.5) ..(1.8,0);
		\draw[blue,postaction={decorate}] (.2,0) ..  controls (-1,-.5) ..(-1.8,0);
		\draw[blue,postaction={decorate}] (.6,0) ..  controls (-.6,-.5) ..(-1.4,0);
		\draw[blue,postaction={decorate}] (1.8,0) ..  controls (1,-.5) ..(-.2,0);
		\fill[black] (-1.8,0) circle (1pt) node[below] {\tiny 1};
		\fill[black] (-1.4,0) circle (1pt) node[below] {\tiny 2};
		\fill[black] (-0.2,0) circle (1pt) node[below] {\tiny n};
		\fill[black] ( .15,0) circle (1pt) node[below] {\tiny n+1};
		\fill[black] ( .6,0) circle (1pt)  node[below]{\tiny n+2};
		\fill[black] ( 1.8,0) circle (1pt) node[below] {\tiny 2n};
		\node[below] at (-.8,0) {$\cdots$};
		\node[below] at (1.2,0) {$\cdots$};
\end{tikzpicture}
	\]

Chord diagrams are permutations of special kind, namely,
involutions without fixed points.
For them, the initial definition coincides with the one above.

Define \emph{Casimir elements} $C_m\in U\gl(N)$, $m=1,2,\dots$:
$$
C_m=w_{\gl(N)}((1,2,\dots,m))=\sum_{i_1,i_2,\dots,i_m=1}^N E_{i_1,i_2}E_{i_2,i_3}\dots E_{i_m,i_1};
$$
associated to the standard cycles $1\mapsto 2\mapsto3\mapsto\dots\mapsto m\mapsto1$.

\begin{theorem}
The center $ZU\gl(N)$ of the universal enveloping algebra $U\gl(N)$ of $\gl(N)$
is identified with the polynomial ring ${\mathbb C}[C_1,\dots,C_N]$.

\end{theorem}

\begin{theorem}[Zhuoke Yang,~\cite{ZY22}]\label{t-gl}
The $w_{\gl(N)}$ invariant of permutations possesses the following properties:
\begin{itemize}
  \item for the empty permutation, the value of $w_{\gl(N)}$ is equal to $1$;
  \item $w_{\gl(N)}$ is multiplicative with respect to concatenation of permutations;
\item {\rm(}\textbf{Recurrence Rule}{\rm)} For the graph of an arbitrary permutation $\alpha$ in~$\S_m$,
and for any two neighboring elements $l,l+1$, of the permuted set $\{1,2,\dots,m\}$, we have
for the values of the $w_{\gl(N)}$ weight system
\begin{equation}\label{e-genrec}
\begin{tikzpicture}[scale=0.8]
\draw (0.6,0) node[v] (V1) {} (1.2,0) node[v] (V2) {};
\draw  (0.5,-.25) node {\scriptsize$\ell$}
 (1.3,-.25) node {\scriptsize$\ell{+}1$};
\draw [->] (0,0) -- (1.8,0);
\draw[->,blue] (0,0.7) .. controls (0.6,0.5) .. (V2);
\draw[->,blue] (V2) .. controls (0.6,0.9) .. (0,1.4);
\draw[->,blue] (1.8,1.4) .. controls (1.2,0.9) .. (V1);
\draw[->,blue] (V1) .. controls (1.2,0.5) .. (1.8,0.7);
\end{tikzpicture}
=
\begin{tikzpicture}[scale=0.8]
\draw  (0.7,0) node[v] (V1) {} (1.1,0) node[v] (V2) {};
\draw [->] (0,0) -- (1.8,0);
\draw[->,blue] (0,0.7) .. controls (0.4,0.4) .. (V1);
\draw[->,blue] (V1) .. controls (0.4,0.9) .. (0,1.4);
\draw[->,blue] (1.8,1.4) .. controls (1.4,0.9) .. (V2);
\draw[->,blue] (V2) .. controls (1.4,0.4) .. (1.8,0.7);
\end{tikzpicture}
+
\begin{tikzpicture}[scale=0.8]
\draw (0.9,0) node[v] (V) {};
\draw [->] (0,0) -- (1.8,0);
\draw[->,blue] (0,0.7) .. controls +(0.5,-0.5) and +(-0.5,-0.5) .. (1.8,0.7);
\draw[->,blue] (1.8,1.4) .. controls +(-0.6,-0.7) .. (V);
\draw[->,blue] (V) .. controls +(-0.3,0.7) .. (0,1.4);
\end{tikzpicture}
 -
\begin{tikzpicture}[scale=0.8]
\draw (0.9,0) node[v] (V) {};
\draw [->] (0,0) -- (1.8,0);
\draw[->,blue] (0,0.7) .. controls (0.4,0.5) .. (V);
\draw[->,blue] (V) .. controls (1.4,0.5) .. (1.8,0.7);
\draw[->,blue] (1.8,1.4) .. controls (0.9,0.2) .. (0,1.4);
\end{tikzpicture}
\end{equation}

\end{itemize}
\end{theorem}

For the special case $\alpha(l+1)=l$, the recurrence looks like follows:
\begin{equation}\label{e-specrec}
	\begin{tikzpicture}[baseline={([yshift=-.5ex]current bounding box.center)},decoration={markings, mark= at position .55 with {\arrow{stealth}}}]
		\draw[->,thick] (-1,0) --  (1,0);
		\fill[black] (-.3,0) circle (1pt) node[below] {\tiny $l$};
		\fill[black] (.3,0) circle (1pt) node[below] {\tiny $l+1$};
		\draw (-.5,.8) node[left] {a};
		\draw (.5,.8) node[right] {b};
		\draw[blue,postaction={decorate}] (-.5,.8) -- (.3,0);
		\draw[blue,postaction={decorate}] (-.3,0) -- (.5,.8);
		\draw[blue,postaction={decorate}] (.3,0) ..controls(0,-.3).. (-.3,0);
	\end{tikzpicture}=
	\begin{tikzpicture}[baseline={([yshift=-.5ex]current bounding box.center)},decoration={markings, mark= at position .55 with {\arrow{stealth}}}]
		\draw[->,thick] (-1,0) --  (1,0);
		\fill[black] (.3,0) circle (1pt) node[below] {\tiny $l+1$};
		\fill[black] (-.3,0) circle (1pt) node[below] {\tiny $l$};
		\draw (-.5,.8) node[left] {a};
		\draw (.5,.8) node[right] {b};
		\draw[blue,postaction={decorate}] (-.5,.8) -- (-.3,0);
		\draw[blue,postaction={decorate}] (.3,0) -- (.5,.8);
		\draw[blue,postaction={decorate}] (-.3,0) ..controls(0,-.3).. (.3,0);
	\end{tikzpicture}+C_1\times
	\begin{tikzpicture}[baseline={([yshift=-.5ex]current bounding box.center)},decoration={markings, mark= at position .55 with {\arrow{stealth}}}]
		\draw[->,thick] (-1,0)  -- (1,0);
		\draw (-.5,.8) node[left] {a};
		\draw (.5,.8) node[right] {b};
		\draw[blue,postaction={decorate}] (-.5,.8) ..controls (0,.4) .. (.5,.8);
	\end{tikzpicture}-N\times
	\begin{tikzpicture}[baseline={([yshift=-.5ex]current bounding box.center)},decoration={markings, mark= at position .55 with {\arrow{stealth}}}]
		\draw[->,thick] (-1,0) --  (1,0);
		\fill[black] (0,0) circle (1pt) node[below] {\tiny $l'$};
		\draw (-.5,.8) node[left] {a};
		\draw (.5,.8) node[right] {b};
		\draw[blue,postaction={decorate}] (-.5,.8) -- (0,0);
		\draw[blue,postaction={decorate}] (0,0) -- (.5,.8);
	\end{tikzpicture}
\end{equation}

Below, we will write the four terms in the recurrence relation for a permutation as $\alpha,\alpha',\beta_1,\beta_2$, so
that the (first) relation looks like 
$$
w_\gl(\alpha)=w_\gl(\alpha')+w_\gl(\beta_1)-w_\gl(\beta_2).
$$
We do not specify the place where the two consecutive elements
are exchanged, but we will do this in each specific case.

\begin{example}\rm
Let us compute the value of $w_\gl$ on the cyclic permutation $(1\ 3\ 2)$ by switching the places of the nodes $2$ and $3$:
\begin{equation*}
  \begin{tikzpicture}[baseline={([yshift=-.5ex]current bounding box.center)},decoration={markings, mark= at position .55 with {\arrow{stealth}}}]
    \draw[->,thick] (-1,0) --  (1,0);
    \fill[black] (-.5,0) circle (1pt) node[below] {\tiny 1};
    \fill[black] (.5,0) circle (1pt) node[below] {\tiny 3};
    \fill[black] (0,0) circle (1pt) node[below] {\tiny 2};
    \draw[blue,postaction={decorate}] (-.5,0) ..controls(0,.3).. (.5,0);
    \draw[blue,postaction={decorate}] (.5,0) ..controls(.25,-.3).. (0,0);
    \draw[blue,postaction={decorate}] (0,0) ..controls(-.25,-.3).. (-.5,0);
    \draw (0,-0.7) node { $(1\ 3\ 2)$};
    \draw (0,0.8) node {};
  \end{tikzpicture}=
  \begin{tikzpicture}[baseline={([yshift=-.5ex]current bounding box.center)},decoration={markings, mark= at position .55 with {\arrow{stealth}}}]
    \draw[->,thick] (-1,0) --  (1,0);
    \fill[black] (-.5,0) circle (1pt) ;
    \fill[black] (.5,0) circle (1pt) ;
    \fill[black] (0,0) circle (1pt) ;
    \draw[blue,postaction={decorate}] (.5,0) ..controls(0,-.3).. (-.5,0);
    \draw[blue,postaction={decorate}] (-.5,0) ..controls(-.25,.3).. (0,0);
    \draw[blue,postaction={decorate}] (0,0) ..controls(.25,.3).. (.5,0);
    \draw (0,-0.7) node { $(1\ 2\ 3)$};
    \draw (0,0.8) node {};
  \end{tikzpicture}+C_1\times
  \begin{tikzpicture}[baseline={([yshift=-.5ex]current bounding box.center)},decoration={markings, mark= at position .55 with {\arrow{stealth}}}]
    \fill[black] (0,0) circle (1pt) ;
    \draw[->,thick] (-1,0)  -- (1,0);
    \draw[blue,postaction={decorate}] (0,0) ..controls (-.5,.7) and (.5,.7)   .. (0,0);
    \draw (0,-0.47) node { $(1)$};
  \end{tikzpicture}-N\times
  \begin{tikzpicture}[baseline={([yshift=-.5ex]current bounding box.center)},decoration={markings, mark= at position .55 with {\arrow{stealth}}}]
    \draw[->,thick] (-1,0) --  (1,0);
    \fill[black] (-.4,0) circle (1pt) node[above] {};
    \fill[black] (.4,0) circle (1pt) node[above] {};
    \draw[blue,postaction={decorate}] (-.4,0) ..controls(0,-.3).. (.4,0);
    \draw[blue,postaction={decorate}] (.4,0) ..controls(0,.3).. (-.4,0);
    \draw (0,0.8) node {};
    \draw (0,-0.7) node { $(1\ 2)$};
  \end{tikzpicture}
\end{equation*}
\begin{eqnarray*}
w_{\gl(N)}((1\ 3\ 2))&=&w_{\gl(N)}((1\ 2\ 3))+C_1\cdot w_{\gl(N)}((1))-N\cdot w_{\gl(N)}((1\ 2))\\
                &=&C_3+C_1^2-NC_2
\end{eqnarray*}
\end{example}

\begin{corollary}
The $\gl(N)$-weight systems, for $N=1,2,\dots$, are combined into
a universal $\gl$-weight system taking values in the ring of
polynomials in infinitely many variables $\C[N;C_1,C_2,\dots]$.

After substituting a given value of~$N$ and an expression of
higher Casimirs $C_{N+1},C_{N+2},\dots$ in terms of the
lower ones $C_1,C_2,\dots,C_N$, this weight system specifies
into the $\gl(N)$-weight system.
\end{corollary}

We say that a function on permutations \emph{is of $\gl$-type}
if it is invariant with respect to the standard shift
and satisfies the same recurrence relation as the $\gl$-weight
system.

\begin{example}
    The number $f(\alpha)$ of faces of a permutation is a function of $\gl$-type. In fact, it is easy to see that
    the number of faces is preserved by replacing the
    first permutation on the left hand side of the recurrence
    with the first one on the right hand side, as well as
    replacing the second permutation on the left with the
    second one on the right.
\end{example}

\section{Rotational Hopf algebras of chord diagrams and permutations}

In this section, we introduce three Hopf algebras, the rotational
Hopf algebra of chord diagrams, 
the rotational Hopf algebra
of permutations, and the rotational Hopf algebra of positive permutations. They have similar structural properties,
and form a chain of 
graded Hopf algebras with respect to inclusion.
Our further goal will be to associate to the $\gl$-weight system
a Hopf algebra homomorphism
of each of these three Hopf algebras to the Hopf algebra of polynomials
in infinitely many variables. 

We start with the rotational Hopf algebra of permutations.
By the \emph{concatenation product of two permutations}
~$\alpha_1\in\S_{m_1}$ and $\alpha_2\in\S_{m_2}$ we mean the
permutation $\alpha_1 \# \alpha_2\in \S_{m_1+m_2}$ whose
directed graph is obtained by attaching the base line of the directed
graph of $\alpha_2$ to the right of that of the directed graph of 
$\alpha_1$, followed by the shift of the numbering of permuted elements:

\[
	\begin{tikzpicture}[baseline={([yshift=-.5ex]current bounding box.center)},decoration={markings, mark= at position .55 with {\arrow{stealth}}}]
		\draw[->,thick] (-2,0)--(0,0);
		\draw[blue,postaction={decorate}] (-1.8,0) ..  controls (-1.4,.3) ..(-1.0,0);
		\draw[blue,postaction={decorate}] (-1.4,0) ..  controls (-1.6,-.3) ..(-1.8,0);
		\draw[blue,postaction={decorate}] (-1.0,0) ..  controls (-.6,.3) ..(-.2,0);
		\draw[blue,postaction={decorate}] (-.6,0) ..  controls (-1,-.3) ..(-1.4,0);
		\draw[blue,postaction={decorate}] (-.2,0) ..  controls (-.4,-.3) ..(-.6,0);

        \fill[black] (-1.8,0) circle (1pt) node[below] {\tiny 1};
		\fill[black] (-1.4,0) circle (1pt) node[below] {\tiny 2};
        \fill[black] (-1.0,0) circle (1pt) node[below] {\tiny 3};
        \fill[black] (-0.6,0) circle (1pt) node[below] {\tiny 4};
		\fill[black] (-0.2,0) circle (1pt) node[below] {\tiny 5};
\end{tikzpicture}
~\#~
	\begin{tikzpicture}[baseline={([yshift=-.5ex]current bounding box.center)},decoration={markings, mark= at position .55 with {\arrow{stealth}}}]
		\draw[->,thick] (0,0)--(1.2,0);
 
		\draw[blue,postaction={decorate}] (.2,0) ..  controls (.6,.3) ..(1,0);
        \draw[blue,postaction={decorate}] (.6,0) ..  controls (.4,-.3) ..(.2,0);
        \draw[blue,postaction={decorate}] (1,0) ..  controls (.8,-.3) ..(.6,0);

		\fill[black] ( .2,0) circle (1pt) node[below] {\tiny 1};
		\fill[black] ( .6,0) circle (1pt)  node[below]{\tiny 2};
		\fill[black] ( 1.0,0) circle (1pt) node[below] {\tiny 3};
\end{tikzpicture}
~=~
	\begin{tikzpicture}[baseline={([yshift=-.5ex]current bounding box.center)},decoration={markings, mark= at position .55 with {\arrow{stealth}}}]
		\draw[->,thick] (-2,0)--(1.2,0);
		\draw[blue,postaction={decorate}] (-1.8,0) ..  controls (-1.4,.3) ..(-1.0,0);
		\draw[blue,postaction={decorate}] (-1.4,0) ..  controls (-1.6,-.3) ..(-1.8,0);
		\draw[blue,postaction={decorate}] (-1.0,0) ..  controls (-.6,.3) ..(-.2,0);
		\draw[blue,postaction={decorate}] (-.6,0) ..  controls (-1,-.3) ..(-1.4,0);
		\draw[blue,postaction={decorate}] (-.2,0) ..  controls (-.4,-.3) ..(-.6,0);
  
		\draw[blue,postaction={decorate}] (.2,0) ..  controls (.6,.3) ..(1,0);
        \draw[blue,postaction={decorate}] (.6,0) ..  controls (.4,-.3) ..(.2,0);
        \draw[blue,postaction={decorate}] (1,0) ..  controls (.8,-.3) ..(.6,0);
        
        \fill[black] (-1.8,0) circle (1pt) node[below] {\tiny 1};
		\fill[black] (-1.4,0) circle (1pt) node[below] {\tiny 2};
        \fill[black] (-1.0,0) circle (1pt) node[below] {\tiny 3};
        \fill[black] (-0.6,0) circle (1pt) node[below] {\tiny 4};
		\fill[black] (-0.2,0) circle (1pt) node[below] {\tiny 5};
		\fill[black] ( .2,0) circle (1pt) node[below] {\tiny 6};
		\fill[black] ( .6,0) circle (1pt)  node[below]{\tiny 7};
		\fill[black] ( 1.0,0) circle (1pt) node[below] {\tiny 8};
\end{tikzpicture}
	\]

Let $\cA_m$ be the vector space (for definiteness, all the vector spaces
we consider are over the ground field $\C$ of complex numbers) freely
spanned by equivalence classes of permutations of~$m$ elements modulo
the following equivalence relation:
\begin{itemize}
    \item any cyclic shift of a permutation is equivalent to it, 
    $\alpha\sim\sigma^{-1}\alpha\sigma$;
    \item the concatenation product of any two permutations is equivalent
    to the concatenation product of any pair of permutations 
    cyclically equivalent to the factors.
\end{itemize}

Denote the equivalence class of a permutation~$\alpha$
by $[\alpha]$. Let $$
\cA=\cA_0\oplus \cA_1\oplus \cA_2\oplus\dots
$$
be the infinite dimensional vector space, which is the direct sum of
the finite dimensional vector spaces $\cA_m$, $m=0,1,2,\dots$,
the vector space~$\cA_m$ being freely spanned by the equivalence
classes of permutations of~$m$ elements. 
The concatenation product of permutations makes $\cA$ into
an infinite dimensional commutative graded algebra.
We introduce the comultiplication $\mu:\cA\to\cA\otimes\cA$ in
the following way:
for a generator $[\alpha]\in \cA_m$ take a permutation $\alpha\in\S_m$
representing this generator and set
$$
\mu([\alpha])=[\mu(\alpha)]=\sum_{I\sqcup J=V(\alpha)}[{\alpha|_{I}]\otimes[\alpha|_{J}}],
$$
where the summation is carried over all the ways to represent the set
$V(\alpha)$ of disjoint cycles of~$\alpha$ into a disjoint
union of two subsets $I,J$, and $\alpha|_{I}$
denotes the restriction of $\alpha$ to the subset 
of elements contained in the cycles from $I$.
Here, following~\cite{KoL23}, we define the \emph{restriction} $\alpha|_U\in\S_{|U|}$
of a permutation~$\alpha\in\S_m$
to a subset $U\subset\{1,2,\dots,m\}$ as a permutation of the elements
of~$U$ preserving their relative cyclic order inside each disjoint 
cycle in~$\alpha$, composed with the subsequent renumbering to $\{1,2,\dots,|U|\}$.

The standard argument shows that the vector space $\cA$ 
endowed with the multiplication $\#$ and the comultiplication $\mu$
is a graded connected commutative cocommutative Hopf algebra.
We call it the \emph{rotational Hopf algebra of permutations}.

\begin{definition}\rm
    An \emph{ascent} of a permutation~$\alpha\in\S_m$
    is a number $i$, $0<i<m$ such that $\alpha(i)>i$. We denote
    the number of ascents in~$\alpha$ by $a(\alpha)$.
    
    A  permutation~$\alpha\in\S_m$ 
    is \emph{positive} if  $a(\alpha)=m-c(\alpha)$ where $c(\alpha)$ is the number of cycles in $\alpha$.
\end{definition}

Obviously, the number of ascents of a permutation is preserved under
its cyclic shifts, whence is well-defined for the equivalence class of a permutation.

The standard cycle $\sigma=(1,2,3,\dots,m)$ is the only positive permutation
consisting of a single cycle. Its number of ascents is $a(\sigma)=m-1$.
Its inverse $\sigma^{-1}=(m,m-1,m-2,\dots,1)$ has the minimal possible 
number of ascents~$1$.
The numbers of ascents of other cyclic permutations are strictly
between these two extremal values.

A permutation is positive if and only if its restriction to each 
disjoint cycle is a standard cycle.
Positive permutations in~$\S_m$  are in one-to-one correspondence 
with partitions of the set $\{1,2,\dots,m\}$: each part of the partition
specifies a single positive disjoint cycle.
In  particular, any  arc diagram,  and, more generally, any involution  is a positive permutation. Remark  that  a cyclic shift of a positive permutation
also is a positive permutation.

Equivalence classes of positive permutations generate
a graded Hopf subalgebra $\cA^+\subset\cA$, 
$$
\cA^+=\cA^+_0\oplus \cA^+_1\oplus \cA^+_2\oplus\dots,$$
$\cA^+_m\subset\cA_m$, $m=0,1,2,\dots,$
which we call
the \emph{rotational Hopf algebra of positive permutations}.

By definition, the Hopf subalgebra
$$
\cB=\cB_0\oplus \cB_1\oplus \cB_2\oplus\dots,
$$
where $\cB_m\subset \cA_m^+ \subset \cA_m$ $m=0,1,2,\dots$ are the vector
subspaces spanned by 
equivalence classes of involutions without fixed points,
is the \emph{rotational Hopf algebra of chord diagrams}.
The definition means, in particular, that all the 
vector spaces~$B_m$ with odd indices~$m$ are zero.

Note that the rotational Hopf algebra of chord diagrams
$\cB$ is different from the conventional Hopf algebra of chord
diagrams~\cite{K93,CDBook12}, since, in contrast to the usual case,
we do not impose the equivalence modulo $4$-term relations.
The latter can also be extended to arbitrary permutations, 
which leads to a corresponding new Hopf algebra structure,
but we do not require it in the present paper.

An equivalence class of a permutation $[\alpha]$ is said to be
\emph{connected} if it cannot be represented as the concatenation product
of two equivalence classes of smaller degree. 
Each of the Hopf algebras $\cA$, $\cA^+$, and $\cB$ is generated by 
connected equivalence classes, and there are no further 
linear relations. This means that the dimension of the 
vector space $P(A_m)$ of primitive elements of degree~$m$ is equal
to the number of connected equivalence classes of permutations in~$\S_m$.

Two more Hopf algebras we will need are the Hopf algebra of polynomials
in infinitely many variables $\C[p_1,p_2,p_3,\dots]$
and the Hopf algebra of graphs~$\cG$.
The comultiplication in 
the Hopf algebra of polynomials is defined on the generators
by $\mu:p_m\mapsto1\otimes p_m+p_m\otimes1$, $m=1,2,3,\dots$. Since the Hopf algebra homomorphisms
we are going to construct will be just filtered, not graded,
we define a filtration rather than grading on this Hopf algebra by setting the degree
of a variable $p_m$ equal to~$m$, $m=1,2,3,\dots$.

The Hopf algebra of graphs~$\cG$ is defined  as an infinite direct sum
$$
\cG=\cG_0\oplus\cG_1\oplus\cG_2\oplus\dots,
$$
where the vector space~$\cG_m$, $m=0,1,2,\dots$, is freely spanned by isomorphism classes of graphs on~$m$ vertices.
The multiplication in it is induced by the disjoint union,
and the comultiplication $\mu:\cG\to\cG\otimes\cG$
is defined on the generators by 
$$
\mu:G\to \sum_{I\sqcup J=V(G)}G|_I\otimes G_J,
$$
where the summation is carried over all (ordered) partitions of
the set $V(G)$ of vertices of~$G$ into two disjoint subsets
and $G_I$ denotes the subgraph of~$G$ induced by the subset
$I\subset V(G)$ of its vertices. Once again, we treat
the Hopf algebra of graphs~$\cG$ as filtered by the
subspaces
$$
\cG_0\subset\cG_0\oplus\cG_1\subset \cG_0\oplus\cG_1\oplus\cG_2\subset\dots
$$
rather than a graded one.

The notion of intersection graph of a chord diagram~\cite{CDL1}
admits a natural extension to arbitrary permutations.
Namely, for a permutation~$\alpha$ of the set~$\{1,2,\dots,m\}$, we say that 
two its disjoint cycles \emph{interlace} if there
are two elements in the first cycle and two elements 
in the second one following the base line in the alternating order.
The \emph{intersection graph} $\gamma(\alpha)$ of a permutation~$\alpha$ is the abstract
graph whose vertices are in one-to-one correspondence with
the disjoint cycles of~$\alpha$, and two vertices 
are connected by an edge iff the corresponding cycles interlace.
Equivalent permutations have isomorphic intersection
graphs. An equivalence class is connected if and only if
its intersection graph is connected.

The mapping $[\alpha]\to\gamma(\alpha)$ 
taking an equivalence class of
permutations to its intersection graph extends by linearity
to a Hopf algebra homomorphism ${\gamma:\cA\to\cG}$ from the 
rotational Hopf 
algebra of permutations to the Hopf algebra of graphs.
Note, however, that this homomorphism is filtered rather
than graded: the grading in the Hopf algebra of graphs is
given by the number of vertices, which corresponds to the number of cycles in a permutation rather than the number of
permuted elements.



Since all the Hopf algebras we consider are generated by their
primitive elements, their algebra homomorphism is a
homomorphism of Hopf algebras iff it takes any primitive
element in the source to a primitive element in the target.

\section{Chromatic substitution}

Our goal in the present section is to define a 
filtered Hopf algebra homomorphism
from the rotational Hopf algebra of permutations~$\cA$
to the Hopf algebra of polynomials in infinitely many variables.
This homomorphism, when specialized to a single variable and 
restricted to the rotational Hopf algebras of positive permutations~$\cA^+$ and of chord diagrams~$\cB$,
provides the chromatic polynomial of the intersection graph
of a positive permutation  or   a chord diagram, respectively.

For a permutation~$\alpha$, we set
\begin{eqnarray*}
X(\alpha)=N^{c(\alpha)-m}w_\gl(\alpha)|_{C_k:=p_kN^{k-1},k=1,2,\dots}.  
\end{eqnarray*}
The value $X(\alpha)$ is a Laurent polynomial in~$N$ and
a polynomial in the variables $p_1,p_2,\dots$, which we call the
\emph{prechromatic invariant}.

The recursion for $w_{\mathfrak{gl}}$, when rewritten in terms of~$X$, acquires the following form:
\begin{itemize}
\item the invariant $X$ takes equal values on all permutations in one equivalence class;
\item it is multiplicative with respect to concatenation product of permutations;
\item the value of $X$ on the standard (positive) $m$-cycle is equal to~$p_m$;
\item if $\ell$ and $\ell+1$ belong to different cycles of~$\alpha$, we have
\begin{equation}\label{eq:rec-X-join}
X\left(~
\begin{tikzpicture}[scale=0.8]
\draw (0.6,0) node[v] (V1) {} (1.2,0) node[v] (V2) {};
\draw  (0.5,-.25) node {\scriptsize$\ell$}
 (1.3,-.25) node {\scriptsize$\ell{+}1$};
\draw [->] (0,0) -- (1.8,0);
\draw[->,blue] (0,0.7) .. controls (0.6,0.5) .. (V2);
\draw[->,blue] (V2) .. controls (0.6,0.9) .. (0,1.4);
\draw[->,blue] (1.8,1.4) .. controls (1.2,0.9) .. (V1);
\draw[->,blue] (V1) .. controls (1.2,0.5) .. (1.8,0.7);
\end{tikzpicture}
 ~\right)=X\left(~
\begin{tikzpicture}[scale=0.8]
\draw  (0.7,0) node[v] (V1) {} (1.1,0) node[v] (V2) {};
\draw [->] (0,0) -- (1.8,0);
\draw[->,blue] (0,0.7) .. controls (0.4,0.4) .. (V1);
\draw[->,blue] (V1) .. controls (0.4,0.9) .. (0,1.4);
\draw[->,blue] (1.8,1.4) .. controls (1.4,0.9) .. (V2);
\draw[->,blue] (V2) .. controls (1.4,0.4) .. (1.8,0.7);
\end{tikzpicture}
 ~\right)+X\left(~
\begin{tikzpicture}[scale=0.8]
\draw (0.9,0) node[v] (V) {};
\draw [->] (0,0) -- (1.8,0);
\draw[->,blue] (0,0.7) .. controls +(0.5,-0.5) and +(-0.5,-0.5) .. (1.8,0.7);
\draw[->,blue] (1.8,1.4) .. controls +(-0.6,-0.7) .. (V);
\draw[->,blue] (V) .. controls +(-0.3,0.7) .. (0,1.4);
\end{tikzpicture}
 ~\right)-X\left(~
\begin{tikzpicture}[scale=0.8]
\draw (0.9,0) node[v] (V) {};
\draw [->] (0,0) -- (1.8,0);
\draw[->,blue] (0,0.7) .. controls (0.4,0.5) .. (V);
\draw[->,blue] (V) .. controls (1.4,0.5) .. (1.8,0.7);
\draw[->,blue] (1.8,1.4) .. controls (0.9,0.2) .. (0,1.4);
\end{tikzpicture}
 ~\right)
\end{equation}
\item if $\ell$ and $\ell+1$ belong to one cycle of~$\alpha$ and $\ell+1\ne\alpha^{\pm1}(\ell)$, we have
\begin{equation}\label{eq:rec-X-cut}
X\left(~
\begin{tikzpicture}[scale=0.8]
\draw (0.6,0) node[v] (V1) {} (1.2,0) node[v] (V2) {};
\draw  (0.5,-.25) node {\scriptsize$k$}
 (1.3,-.25) node {\scriptsize$k{+}1$};
\draw [->] (0,0) -- (1.8,0);
\draw[->,blue] (0,0.7) .. controls (0.6,0.5) .. (V2);
\draw[->,blue] (V2) .. controls (0.6,0.9) .. (0,1.4);
\draw[->,blue] (1.8,1.4) .. controls (1.2,0.9) .. (V1);
\draw[->,blue] (V1) .. controls (1.2,0.5) .. (1.8,0.7);
\end{tikzpicture}
 ~\right)=X\left(~
\begin{tikzpicture}[scale=0.8]
\draw  (0.7,0) node[v] (V1) {} (1.1,0) node[v] (V2) {};
\draw [->] (0,0) -- (1.8,0);
\draw[->,blue] (0,0.7) .. controls (0.4,0.4) .. (V1);
\draw[->,blue] (V1) .. controls (0.4,0.9) .. (0,1.4);
\draw[->,blue] (1.8,1.4) .. controls (1.4,0.9) .. (V2);
\draw[->,blue] (V2) .. controls (1.4,0.4) .. (1.8,0.7);
\end{tikzpicture}
 ~\right)+N^{-2}\cdot\left(X\left(~
\begin{tikzpicture}[scale=0.8]
\draw (0.9,0) node[v] (V) {};
\draw [->] (0,0) -- (1.8,0);
\draw[->,blue] (0,0.7) .. controls +(0.5,-0.5) and +(-0.5,-0.5) .. (1.8,0.7);
\draw[->,blue] (1.8,1.4) .. controls +(-0.6,-0.7) .. (V);
\draw[->,blue] (V) .. controls +(-0.3,0.7) .. (0,1.4);
\end{tikzpicture}
 ~\right)-X\left(~
\begin{tikzpicture}[scale=0.8]
\draw (0.9,0) node[v] (V) {};
\draw [->] (0,0) -- (1.8,0);
\draw[->,blue] (0,0.7) .. controls (0.4,0.5) .. (V);
\draw[->,blue] (V) .. controls (1.4,0.5) .. (1.8,0.7);
\draw[->,blue] (1.8,1.4) .. controls (0.9,0.2) .. (0,1.4);
\end{tikzpicture}
 ~\right)\right)
\end{equation}
\item in the case $\alpha(\ell+1)=\ell$, we have
\begin{equation}\label{eq:rec-X-cut0}
X\left(~
\begin{tikzpicture}[scale=0.8]
\draw (0.6,0) node[v] (V1) {} (1.2,0) node[v] (V2) {};
\draw  (0.5,-.25) node {\scriptsize$k$}
 (1.3,-.25) node {\scriptsize$k{+}1$};
\draw [->] (0,0) -- (1.8,0);
\draw[->,blue] (0,1) .. controls (0.7,0.6) .. (V2);
\draw[->,blue] (V1) .. controls (1.1,0.6) .. (1.8,1);
\draw[->,blue] (V2) .. controls (0.9,0.2) .. (V1);
\end{tikzpicture}
 ~\right)= X\left(~
\begin{tikzpicture}[scale=0.8]
\draw  (0.6,0) node[v] (V1) {} (1.2,0) node[v] (V2) {};
\draw [->] (0,0) -- (1.8,0);
\draw[->,blue] (0,1) .. controls (0.4,0.6) .. (V1);
\draw[->,blue] (V2) .. controls (1.4,0.6) .. (1.8,1);
\draw[->,blue] (V1) .. controls (0.9,0.3) .. (V2);
\end{tikzpicture}
 ~\right)+
   p_1N^{-2}\cdot X\left(~
\begin{tikzpicture}[scale=0.8]
\draw [->] (0,0) -- (1.8,0);
\draw[->,blue] (0,1) .. controls (0.9,0) .. (1.8,1);
\end{tikzpicture}
 ~\right)-
X\left(~
\begin{tikzpicture}[scale=0.8]
\draw (0.9,0) node[v] (V) {};
\draw [->] (0,0) -- (1.8,0);
\draw[->,blue] (V) .. controls (1.2,0.5) .. (1.8,1);
\draw[->,blue] (0,1) .. controls (0.6,0.5) .. (V);
\end{tikzpicture}
 ~\right)
\end{equation}
\end{itemize}

As an immediate corollary of these relations we obtain

\begin{lemma}
The prechromatic invariant takes values in the ring of polynomials in $N^{-2},p_1,p_2,\dots$. In other words, its value on any permutation~$\alpha$ can be written as a finite sum
$$
X(\alpha)=X_0(\alpha)+N^{-2}X_1(\alpha)+N^{-4}X_2(\alpha)+\dots
$$ 
where each $X_k(\alpha)$ is a polynomial in the variables
$p_1,p_2,\dots$.
\end{lemma}

We are mostly interested in the leading term $X_0$ of this expansion.
The following assertion is the main result of the present paper.


\begin{theorem}\label{th-HAh}
    The value $X_0(\alpha)$ on a permutation 
    $\alpha\in\S_m$ 
    extends by linearity to a
    filtered Hopf algebra homomorphism $X_0:\cA\to\C[p_1,p_2,\dots]$.
\end{theorem}

\begin{corollary}\label{c-cp}
    The restriction of $X_0$ to the rotational Hopf algebras
    of positive  permutations  $\cA^+\subset\cA$  and  of chord
    diagrams $\cB\subset\cA$ is a filtered Hopf algebra homomorphism.
    Under the substitution $p_m:=x$, $m=1,2,\dots$ its value
    on an equivalence class of a positive permutation
    or a  chord diagram becomes
    the chromatic polynomial of the corresponding intersection graph.
  \end{corollary}

We   call the substitution $p_m:=x $ to~$X$
the  \emph{chromatic substitution}.

\begin{example}
The positive permutation $(1,3,5,8,10)(2,4,7)(6,9,11)$  consists  of
three disjoint cycles. Any two of these cycles interlace,
whence its  intersection graph is the
complete graph~$K_3$. The  chromatic substitution on this permutation is
\begin{eqnarray*}
&&X((1,3,5,8,10)(2,4,7)(6,9,11))(N;x,x,\dots)\\
&=&x(x-1)(x-2)\\
&+&x(x-1) (3x-5) N^{-2}\\
&+&2x(x-1)(3x-2) N^{-4}\\
&+&x(x-1)(5x-1) N^{-6}\\
&+&x^2(2x-1)N^{-8}.
\end{eqnarray*}
The coefficient of the leading term in~$N$ is the chromatic polynomial $\chi_{K_3}(x)$.
\end{example}

\begin{example}
For the value of the $\gl$-weight system on the chord diagram
of order~$5$ whose intersection
graph is the complete graph~$K_5$, we have
\begin{eqnarray*}
w_\gl(K_5)&=&{\scriptstyle  24C_2N^4+(24C_3-50C_2^2-24C_1^2)N^3}\\
&&{\scriptstyle  -(24C_4+10C_2C_3-35C_2^3-70C_1^2C_2+72C_1C_2-32C_2)N^2}\\
&&{\scriptstyle  +(10C_2C_4+96C_1C_3-10C_2^4-50C_1^2C_2^2+30C_1C_2^2-82C_2^2-20C_1^4+48C_1^3-32C_1^2)N}\\
&&{\scriptstyle -40C_1C_2C_3+C_2^5+10C_1^2C_2^3+30C_2^3+15C_1^4C_2-20C_1^3C_2+10C_1^2C_2},
\end{eqnarray*}
which, under the chromatic substitution, becomes
\begin{eqnarray*}
X(K_5)(N;x,x,\dots)
&=&x(x-1)(x-2)(x-3)(x-4)\\
&+&2x (x-1) (5x^3-20x^2+25x-16)N^{-2}\\
&+&x^2 (15 x^3-40x^2+58x-32)N^{-4}.
\end{eqnarray*}

\end{example}

It would be interesting to know the combinatorial meaning
of the coefficients of the chromatic substitution at
negative powers of~$N$. In particular, whether they
depend on the intersection graph of a chord diagram (or
a positive permutation) only. 

\begin{theorem}\label{th-csde}
    The value of the chromatic substitution on a nonpositive 
    permutation~$\alpha\in\S_m$ is a polynomial  in~$N^{-2}$
    with zero free term.
\end{theorem}

\begin{example}
    For the inverse standard cycle $(132)$, the chromatic substitution
    to the value of the $\gl$-weight system yields
    $$X_((132))(N;x,x,\dots)=N^{-2}(C_3+C_1^2-NC_2)|_{C_k:=xN^{k-1}}=x^2N^{-2}.$$
    Its free term is~$0$.
\end{example}

The following further 
specializations of the chromatic substitution are useful.
It is shown in~\cite{KoL23} that for the substitution $C_k:=N^{k-1}$,
$k=1,2,\dots$, the universal $\gl$-weight system takes a permutation~$\alpha$ to $N^{f(\alpha)-1}$,
where $f(\alpha)$ is the number of faces of~$\alpha$, $w_\gl(\alpha)|_{C_k:=N^{k-1}}=N^{f(\alpha)-1}$.

\begin{lemma} We have the following specializations of the universal
$\gl$-weight system:
\begin{itemize}
\item   for the substitution~$C_k:=N^{k+1}$, $k=1,2,\dots$,
    the universal $\gl$-weight system takes a permutation~$\alpha\in\S_m$ to the 
    monomial $N^{m+c(\alpha)}$, $w_\gl(\alpha)|_{C_k:=N^{k+1}}=N^{m+c(\alpha)}$;

\item for the substitution~$N:=1,C_k:=x$, $k=1,2,\dots$,
    the  universal $\gl$-weight system takes a permutation~$\alpha\in\S_m$ to the 
    monomial $x^{c(\alpha)}$, $w_\gl(\alpha)|_{N:=1,C_k:=x}=x^{c(\alpha)}$. 
    
\end{itemize}
Recall that  $c(\alpha)$ denotes the number of cycles in~$\alpha$.
    
\end{lemma}

Indeed, the first assertion of the lemma is true for the standard cycles. Moreover, the
function $\alpha\mapsto N^{m+c(\alpha)}$ is a function of $\gl$-type, that is,
it satisfies the $\gl$ recurrence relation~(\ref{e-genrec}): the value $m+c(\alpha)$ is the same
for both~$\alpha$ and~$\alpha'$, as well as for both $\beta_1$ and $\beta_2$.

For the second assertion, the argument is the same.

\section{Proof of the main results}

In this section we give a proof of the main results of the paper.
We start with defining a linear order on the set of 
equivalence classes of permutations in~$\S_m$, that is,
on the set of basic vectors
in the spaces $\cA_m$. Then we argue by a double induction
with respect to~$m$ and this linear ordering.

\subsection{A linear ordering on the set of equivalence classes
of permutations}

Associate to the equivalence class $[\alpha]$ 
of a permutation~$\alpha$ its
canonical form in the following
way. Write each permutation in the equivalence class
as a product of disjoint cycles, where for each cycle 
we chose the lexicographically minimal presentation (the one with the smallest first element),
and write the cycles in the lexicographic order (that is, the order in which their
first elements increase). Now, among all the permutations 
in the equivalence class, choose the permutation with the lexicographically smallest canonical
form. This will be the \emph{canonical form}~$\bar\alpha$ of~$[\alpha]$. The lexicographic ordering of canonical
forms defines a linear ordering on the basic vectors in~$A_m$.

\subsection{Proof of Theorem~\ref{th-HAh}}

 We argue by induction over~$m$, the number of permuted elements,
 and then by the linear ordering of the basic elements of~$\cA_m$
 defined above. 

Suppose we have proved already the assertion for all
permutations of less than~$m$
elements, so that we know that for each permutation~$\alpha$ of at most~$m-1$
elements
the leading term $X_0(\alpha)$ in~$N$ is a linear polynomial 
    in~$p_1,p_2,\dots,p_{i}$.


The comultiplicativity of
the mapping~$X_0:\cA\to\C[p_1,p_2,\dots]$ can be reformulated in
the following way. Consider the mapping
$$
(X_0\otimes X_0)\circ\mu:\cA\to \C[p_1,p_2,\dots]\otimes\C[p_1,p_2,\dots],
$$
which we denote by~$\Xi$. For a given permutation $\alpha$,
it looks like
$$
\Xi:\alpha\mapsto \sum_{I\sqcup J=V(\alpha)}
X_0(\alpha|_I)\otimes X_0(\alpha|_J),
$$
where the summation is carried over all ordered partitions of the
set~$V(\alpha)$ of disjoint cycles of~$\alpha$ into two
disjoint subsets. 

\begin{lemma}\label{l:Xi}
The mapping~$\Xi$ behaves as follows with respect to the recurrence
relations applied to a permutation~$\alpha$:
\begin{itemize}
    \item under~(\ref{eq:rec-X-join}), with $X_0(\alpha)=X_0(\alpha')+
    X_0(\beta_1)-X_0(\beta_2)$, we have
    $$\Xi(\alpha)=\Xi(\alpha')+
    \Xi(\beta_1)-\Xi(\beta_2);$$
    \item under~(\ref{eq:rec-X-cut}), with $X_0(\alpha)=X_0(\alpha')+
   N^{-2}( X_0(\beta_1)-X_0(\beta_2))$, we have
    $$\Xi(\alpha)=\Xi(\alpha');$$
        \item under~(\ref{eq:rec-X-cut0}), with $X_0(\alpha)=X_0(\alpha')+
    N^{-2}X_0(\beta_1)-X_0(\beta_2)$, we have
    $$\Xi(\alpha)=\Xi(\alpha')-\Xi(\beta_2).$$
\end{itemize}

\end{lemma}

Before proving Lemma~\ref{l:Xi} let us deduce Theorem~\ref{th-HAh}
from it.

The comultiplicativity of~$X_0$ is
equivalent to the identity
$$
\Xi(\alpha)=X_0(\alpha)|_{p_k:=p_k\otimes 1+1\otimes p_k}.
$$
Since $X_0(C_k)=p_k$, $k=1,2,3,\dots$, the comultiplicativity property is 
equivalent to the assertion that the image of the mapping
$$
\Xi:\cA\to \C[p_1,p_2,\dots]\otimes\C[p_1,p_2,\dots]
$$
coincides with the subring generated by the elements $p_k\otimes 1+1\otimes p_k$.

Suppose $m$ is the minimal number for which
the assertion of Theorem~\ref{th-HAh} is not proved yet. Let
$[\alpha]$ be the lexicographically minimal basic vector
in~$A_m$ possessing this property, and let~$\bar\alpha$
be the corresponding canonical form.
Obviously, $\bar\alpha$
is connected and contains at least~$2$ disjoint cycles. 
Take the minimal number~$k$, $1\le k<m$ such that
$\bar\alpha(k)>k+1$. Applying the recurrence relation to the
pair $\bar\alpha(k)-1,\bar\alpha(k)$ of permuted elements, we fall
into one of the three cases of Lemma~\ref{l:Xi}.
Since in each of the cases each of the permutations 
$\bar\alpha',\beta_1,\beta_2$ are smaller than $\bar\alpha$,
the induction argument proves the Theorem.

{\bf Proof of Lemma~\ref{l:Xi}.}

\begin{itemize}
    \item If the pair $\ell-1,\ell$ is like in recursion~(\ref{eq:rec-X-join}), so that $\ell-1$ and $\ell$ belong to two different
    disjoint cycles of~$\alpha$, then let $v_1\in V(\alpha)$
    be the cycle containing~$\ell$, and $v_2\in V(\alpha)$ 
    be the cycle containing~$\ell-1$.  The disjoint cycles~$V(\alpha')$
    naturally correspond one-to-one to the disjoint cycles
    $V(\alpha)$, while each of the permutations $\beta_1,\beta_2$
    has one disjoint cycle less: the pair of cycles $v_1,v_2$ is replaced by a single cycle. Now, we have
\begin{eqnarray*}
   && \Xi(\alpha)-\Xi(\alpha')-\Xi(\beta_1)+\Xi(\beta_2)=\\
   && \sum_{{I\sqcup J=V(\alpha)\atop \{v_1,v_2\}\subset I}}
    X_0(\alpha|_I)\otimes(X_0(\alpha|_J)-X_0(\alpha'|_J)
    -X_0(\beta_1|_J)+X_0(\beta_2|_J))\\
   && \sum_{{I\sqcup J=V(\alpha)\atop \{v_1,v_2\}\subset J}}
   (X_0(\alpha|_I)-X_0(\alpha'|_I)
    -X_0(\beta_1|_I)+X_0(\beta_2|_I)) \otimes X_0(\alpha|_J)\\
    &&\sum_{{I\sqcup J=V(\alpha)\atop v_1\in I,v_2\in J
    \text{ or }v_2\in I,v_1\in J}}(X_0(\alpha|_I)\otimes
    X_0(\alpha|_J)-X_0(\alpha'|_I)
    \otimes X_0(\alpha'|_J)).
\end{eqnarray*}
Each of the three summands on the right is zero: for the first two
sums the expression in brackets is~$0$ due to recurrence~(\ref{eq:rec-X-join}) applied to the permutations $\alpha|_J$ and $\alpha|_I$,
respectively,
while for the third sum, it vanishes because $\alpha|_I=\alpha'|_I,
\alpha|_J=\alpha'|_J$.
 \item If the pair $\ell-1,\ell$ is like in recurrence~(\ref{eq:rec-X-cut}), so that $\ell-1$ and $\ell$ belong to the same
    disjoint cycle of~$\alpha$, then let $v\in V(\alpha)$
    be the cycle containing this pair of elements.  
    The disjoint cycles~$V(\alpha')$
    naturally correspond one-to-one to the disjoint cycles
    $V(\alpha)$. We claim that for an arbitrary subset $I\subset V(\alpha)$ the values $X_0(\alpha|_I)$ and $X_0(\alpha'|_I)$
    coincide. Indeed, if $I$ does not contain~$v$, then 
    $\alpha|_I=\alpha'|_I$. If~$I$ contains~$v$, then the 
    two restrictions coincide no longer. However, in this case
    the permutations $\alpha|_I$ and $\alpha'|_I$ are related by 
    the same recurrence~(\ref{eq:rec-X-cut}), and the equality
    follows from the induction hypothesis.
     \item If the pair $\ell-1,\ell$ is like in recursion~(\ref{eq:rec-X-cut0}), so that $\ell-1$ and $\ell$ are two successive elements
     in the same disjoint cycle, then let $v\in V(\alpha)$
    be this cycle.  The disjoint cycles~$V(\alpha')$ and $V(\beta_2)$
    naturally correspond one-to-one to the disjoint cycles
    $V(\alpha)$, while the permutation $\beta_1$
    has one disjoint cycle more: the cycle $v$ is replaced by two cycles. Thus, for any subset $I\subset V(\alpha)$ not containing~$v$, we have $\alpha|_I=\alpha'|_I=\beta_2|_I$.
    As a consequence,
\begin{eqnarray*}
   && \Xi(\alpha)-\Xi(\alpha')+\Xi(\beta_2)=\\
   && \sum_{{I\sqcup J=V(\alpha)\atop v\in I}}
   (X_0(\alpha|_I)-X_0(\alpha'|_I)
   +X_0(\beta_2|_I))\otimes  X_0(\alpha|_J)\\
   && \sum_{{I\sqcup J=V(\alpha)\atop v\in J}}
    X_0(\alpha|_J)\otimes(X_0(\alpha|_J)-X_0(\alpha'|_J)
   +X_0(\beta_2|_J)).
\end{eqnarray*}
Each of the two summands on the right vanishes: the expression in brackets is~$0$ due to recurrence~(\ref{eq:rec-X-cut0}) applied to the permutations $\alpha|_J$ and $\alpha|_I$.
\end{itemize}
Lemma~\ref{l:Xi} is proved, which completes the proof of Theorem~\ref{th-HAh}.

\subsection{Values of~$X_0$ on cyclic permutations}

In this section we compute explicitly the value of~$X_0$
on cyclic permutations.

\begin{proposition}
    If $\alpha\in\S_m$ is a cyclic permutation, that is, $c(\alpha)=1$,
    then the value $X_\alpha$ of the prechromatic invariant on it
    has the form 
    $$
    X(\alpha)=X_0(p_1,p_2,\dots)+
    X_1(p_1,p_2,\dots)N^{-2}+\dots,
    $$
where $X_0$ is a linear polynomial in $p_1,\dots,p_m$,
$X_1$ is at most a quadratic polynomial in $p_1,\dots,p_{m-2}$,
$X_2$ is at most a cubic polynomial in~$p_1,\dots,p_{m-4}$, \dots.
      For a cyclic permutation $\alpha\in\S_m$ having $m-k$ ascents, $a(\alpha)=m-k$,
    the coefficient $X_0$ is 
    $$X_0(\alpha)=p_{m}-{k-1\choose1}p_{m-1}+{k-1\choose2}p_{m-2}-\dots+(-1)^{k-1}{k-1\choose k-1}p_{m-k+1}.$$
\end{proposition}

\begin{remark}
 In particular,  for $m>1$ the coefficient of~$p_1$ in the linear polynomial  $L_1$  is~$0$, so that it depends  only  on
$p_2,\dots,p_m$; for  $k>1$, under   the  substitution
$p_k=x$ the polynomial $X_0(\alpha)$ vanishes.   
\end{remark}

Let $[\alpha]$ be the minimal cyclic permutation for which
the assertion of the Proposition is not proved yet. Take its
canonical form~$\bar\alpha\in\S_m$ and let~$\ell$, $1\le\ell< m-1$ be the minimal element such that $\bar\alpha(\ell)>\ell+1$.
(If such an element~$\ell$ does not exist, then~$\bar\alpha$ is
a standard cycle, for which the assertion of the Proposition is obvious.)
There are two possibilities:
\begin{itemize}
    \item either $\bar\alpha(\bar\alpha(\ell))=\bar\alpha(\ell)-1$;   
    \item or $\bar\alpha(\bar\alpha(\ell))\ne\bar\alpha(\ell)-1$.
\end{itemize}

In the first case, recurrence~(\ref{eq:rec-X-cut0})
is applicable 
to the points $\bar\alpha(\ell)-1,\bar\alpha(\ell)$.
The number of ascents in the permutation $\bar\alpha'$,
the first permutation 
on the right, is one more than $a(\bar\alpha)$, and, by the induction hypothesis, the polynomial $X_0({\bar\alpha'})$ for it is
$$p_m-{m-a(\bar\alpha)-2\choose1}p_{m-1}+{m-a(\bar\alpha)-2\choose2}p_{m-2}-\dots.$$
The second permutation $\beta_1$ on the right-hand side
is a permutation of
$m-2$ elements; by the induction hypothesis, after multiplication
by~$p_1N^{-2}$ the value of~$X$ on it does not
contribute to the coefficient of $N^{0}$. The second 
element $\beta_2$ is a cyclic permutation of~$m-1$ elements, with
$a(\bar\alpha)$ ascents. By the induction hypothesis,  
the value of~$X_0$ on it is
$$p_{m-1}-{m-a(\bar\alpha)-2\choose1}p_{m-2}+{m-a(\bar\alpha)-2\choose2}p_{m-3}-\dots$$
to the coefficient of $N^{m-1}$. After subtracting this
value from the contribution of~$\bar\alpha'$,
we deduce that the value~$X_0(\bar\alpha)$ is as required.

In the second case, apply recurrence~(\ref{eq:rec-X-cut})
to the points $\bar\alpha(\ell)-1,\bar\alpha(\ell)$.
Then $\bar\alpha'$, the first term on the right, is a cyclic permutation
of $m$ elements whose number of ascents coincides
with $a(\bar\alpha)$. The two other
permutations on the right $\beta_1,\beta_2$ 
are permutations of $m-1$ elements, each consisting of two
disjoint cycles. Therefore, the polynomial~$X$
on none of them contributes to free term in $N$, hence the
assertion of the lemma is true for $[\alpha]$ as well.

The property that $X_k$ is a polynomial of degree at most~$k+1$
in~$p_i$ also 
follows immediately from the induction argument above.
The Proposition is proved.

\subsection{Proof  of  Corollary~\ref{c-cp}}

The proof of the fact that the leading term in~$N$ of
the chromatic substitution on
positive permutations coincides with the chromatic polynomial of
their intersection graphs splits into two steps.
First, we need to prove that this leading term depends on
the intersection graph only. We show this by proving that it satisfies
the deletion-contraction relation.
In addition, we prove that for a nonpositive permutation $\alpha$,
the free term in the chromatic
substitution is~$0$.

The proof  proceeds by induction on the number~$m$
of permuted elements, and the ordering on the equivalence classes
$[\alpha]$ of permutations defined above. Suppose we have already
proved Corollary~\ref{c-cp} for all positive permutations of
at most~$m-1$ elements. For the cyclic positive 
permutation~$\sigma\in\S_m$, the assertions are also valid.

Let $[\alpha], \alpha\in\S_m$ be the minimal equivalence class of a 
positive permutation, for which Corollary~\ref{c-cp} is not proved yet.
It contains at least two disjoint cycles.
Take its canonical representative $\bar\alpha$, let~$v_1$ 
be the disjoint cycle containing~$1$, and let $\ell$,
$2\le \ell\le m-1$ be the minimal element of the cycle~$v_1$
such that $\ell-1$ does not belong to this cycle.
Let us denote the cycle of~$\bar\alpha$ which contains~$\ell-1$
by~$v_2$.
Apply recurrence~(\ref{e-genrec}) to the pair of points $\ell-1,\ell$
and then apply the chromatic substitution to all its four terms.
There are two possibilities:
\begin{itemize}
    \item the intersection graph $\gamma(\bar\alpha)$ of the permutation $\bar\alpha$ and the intersection graph $\gamma(\bar\alpha')$ of the permutation~$\bar\alpha'$ are isomorphic to one another;
    \item these two intersection graphs are distinct.
\end{itemize}
In the second case, the intersection graph $\gamma(\bar\alpha')$
is obtained from $\gamma(\bar\alpha)$ by either deleting or
adding the edge connecting~$v_1$ and $v_2$.

Both permutations $\beta_1,\beta_2$ on the right-hand side of
recurrence~(\ref{e-genrec}) are permutations of~$m-1$ elements
containing $c(\alpha)-1$ disjoint cycles. If the intersection graphs
$\gamma(\bar\alpha)$ and $\gamma(\bar\alpha')$ are isomorphic,
then, by the induction hypothesis, the contribution of the chromatic substitution to the coefficient of $N^{0}$
on both~$\beta_1$ and~$\beta_2$ is~$0$. Indeed, each of these two
permutations contains a nonpositive disjoint cycle, the one
obtained by merging the cycles $v_1,v_2$ of~$\bar\alpha$.

Now suppose the graphs $\gamma(\bar\alpha)$ and $\gamma(\bar\alpha')$
are not isomorphic.  Then if $v_1,v_2$ are connected in $\gamma(\bar\alpha)$ and not in $\gamma(\bar\alpha')$, then the permutation
$\beta_1$ is not positive, while $\beta_2$ is. In this case the 
contribution of the chromatic substitution on~$\beta_1$ to
the coefficient of~$N^{0}$ is~$0$, while
the intersection graph $\gamma(\beta_2)$ of the permutation $\beta_2$ is the result
of contracting the edge $(v_1,v_2)$ in the intersection graph $\gamma(\bar\alpha)$. The formula for the coefficient of $N^{0}$ in the chromatic susbstitution for~$\bar\alpha$ then
becomes the contraction-deletion relation for the chromatic polynomial,
which proves that this coefficient is, indeed, the chromatic
polynomial of $\gamma(\bar\alpha)$.
If $v_1,v_2$ are not connected in $\gamma(\bar\alpha)$
and are connected in $\gamma(\bar\alpha')$, the argument is the same,
with the exchange of the roles of the permutations $\beta_1,\beta_2$.

In the case the permutation $\bar\alpha$ is not positive,
all the three permutations $\bar\alpha',\beta_1,\beta_2$
are not positive as well. Hence, by the induction hypothesis,
none of the three permutations
contributes to the coefficient of $N^{m-c(\alpha)}$
in the chromatic substitution, whence
the same is true for the permutation~$\bar\alpha$.
This completes the proof of Corollary~\ref{c-cp}.

It is worth to recall that
the chromatic polynomial can be considered as a filtered Hopf algebra homomorphism
$\chi:\cG\to\C[x]$, see, e.g.~\cite{KL22}. In our current
setting, it can be considered also as a filtered Hopf algebra
homomorphism $\chi:\cA^+\to\C[x]$ from the rotational Hopf algebra
of positive permutations to the Hopf algebra of polynomials in
a single variable.

\section{Further results}


In this section, we present several results about the chromatic substitution, which
both clarify its nature and simplify its computation.

\begin{theorem}\label{th-cplc}
If $\alpha$ is a nonpositive long cycle of length~$m$, then the
degree in~$N$ of the polynomial $X(\alpha)(N;x,x,\dots)$ is less than~$-1$.
In particular, for cyclic permutations with $m-2$ ascents, this polynomial
has degree $-2$ in~$N$.
\end{theorem}

\begin{remark}
   The chromatic substitution on cyclic permutations with less
   than~$m-2$ ascents also can provide polynomials of degree~$-2$ in~$N$.
   For example, the permutation $(1,2, 3, 6,5,4)\in\S_6$ 
   has~$m-3=3$ ascents, while the chromatic substitution on it
   yields
   $$
   X({(1,2, 3, 6,5,4)})(N;x,x,\dots)=x^2N^{-2}.
   $$
\end{remark}

\begin{example}
The long cycle $\alpha=(1,3,5,2,4)$  is nonpositive;
the number of faces in it is $f(\alpha)=1$. The  chromatic substitution on this permutation is
\begin{eqnarray*}
X_\alpha(N;x,x,\dots)
&=&x(x-1) N^{-2}+x^2N^{-4}.
\end{eqnarray*}

\end{example}

\begin{lemma}
    The value of the chromatic substitution on a cyclic permutation in~$\S_m$
    with $m-2$ ascents is $x^2N^{m-3}$.
\end{lemma}

{\bf Proof.} We proceed by induction on~$m$. Suppose the assertion of
the lemma is true for all cyclic permutations of up to $m-1$ elements.
Each cyclic permutation of positivity~$m-4$ can be obtained
from the positive cyclic permutation by a composition of exchanging  
the pairs of elements $(m-1,m)$, $(m-2,m-1)$, \dots. Applying
the recurrence rule for $\gl$ to the first
of this exchanges and using the induction assumption, we conclude that
the value of the chromatic substitution on the cyclic permutation
$(1,2,\dots,m-2,,m,m-1)$ indeed is $x^2N^{m-3}$. Further exchanges
preserve this value since for them the right hand side of the 
recurrence relation vanishes. $\square$

\begin{theorem}
       All the coefficients but one in $X(\alpha)(N;x,x,\dots)$,
        for an arbitrary permutation~$\alpha$ of at least one element, 
        are divisible by $x(x-1)$.
    The only exception is the coefficient of $N^{f(\alpha)-m}$;
    it is divisible by~$x$ and its value
   at $x=1$ is~$1$.    
\end{theorem}

Indeed, we know that the composition of the chromatic substitution with $x:=1$
makes its value on a permutation~$\alpha$ into~$N^{f(\alpha)-m}$, while substitution
$x:=0$ makes it into~$0$ for all permutations but the empty one.

\begin{theorem}
        Let~$\alpha\in \S_m$ be a permutation such that $\alpha(1)=2$. Then 
        $$X(\alpha)(N;x,x,\dots)=X({\alpha|_{\{2,3,\dots,m\}}})(N;x,x,\dots).$$
\end{theorem}

Note that since permutations are considered up to cyclic shifts, 
a similar assertion is valid for any permutation $\alpha\in\S_m$
such that $\alpha(k)\equiv k+1\mod~m$ for some~$k\in\{1,\dots,m\}$.
This means that if a permutation takes an element to the next one in the cyclic order,
then we can replace this pair of consecutive elements by a single one,
the value of the chromatic substitution being preserved.
The proof repeats that for the main theorem for permutations.

Associate to a permutation $\alpha$ the \emph{weighted intersection graph}
$\gamma^w(\alpha)$, which is a simple graph with ${\mathbb N}$-weighted
vertices, in the following way: the underlying graph of $\gamma^w(\alpha)$
coincides with $\gamma(\alpha)$, and the weight of a vertex is the length of
the cycle in~$\alpha$, which corresponds to this cycle.


        

\end{document}